\documentclass[10pt,a4paper]{amsart}

\usepackage{amsfonts,amssymb,latexsym}

\def\LaTeX{L\kern-.36em\raise.3ex\hbox{a}\kern-.15em
    T\kern-.1667em\lower.7ex\hbox{E}\kern-.125emX}

\def\dst{\displaystyle}

\newcommand{\Be}{\begin{equation}}
\newcommand{\Ee}{\end{equation}}

\def\udots{\raisebox{0.35mm}{\mbox{.}\raisebox{1.1mm}{\mbox{\hskip0.2mm.}}\raise
box{2.2mm}{\mbox{\hskip0.3mm.}}}}

\def\C{{\mathbb{C}}}

\def\Q{{\mathbb{Q}}}
\def\R{{\mathbb{R}}}

\def\Z{{\mathbb{Z}}}

\newcommand{\ent}[1]{{\left[{#1}\right]}}
\newcommand{\abs}[1]{{\left|{#1}\right|}}
\newcommand{\scal}[1]{{\left\langle{#1}\right\rangle}}

\newenvironment{notation}[1][]{\vskip1pt\noindent\rm\textit{Notation}:\
 }{\rm\vskip1pt}
\newenvironment{remark}[1][]{\vskip1pt\noindent\rm\textit{Remark}:\
}{\rm\vskip1pt}

\newtheorem{lemma}{Lemma}[section]
\newtheorem{proposition}[lemma]{Proposition}
\newtheorem{theorem}[lemma]{Theorem}

\begin{document}
\title[Distances sets that are a shift of the integers]{Distances sets that are a shift of the integers
and Fourier basis for planar convex sets}
\author{Alex IOSEVICH \& Philippe JAMING}
\address{A.~I.~: Department of Mathematics\\ University of Missouri-Columbia\\
Columbia\\ Missouri 65211\\ USA}
\email{iosevich@math.missouri.edu}
\address{P.~J.\,: MAPMO-F\'ed\'eration Denis Poisson\\ Universit\'e d'Orl\'eans\\ BP 6759\\ F 45067 ORLEANS Cedex 2\\
FRANCE}
\email{Philippe.Jaming@univ-orleans.fr}


\thanks{Research partially financed by : {\it European Commission}
Harmonic Analysis and Related Problems 2002-2006 IHP Network
(Contract Number: HPRN-CT-2001-00273 - HARP) and a grant from the National Science Foundation}

\begin{abstract}
The aim of this paper is to prove that if a planar set $A$ has a difference set $\Delta(A)$
satisfying $\Delta(A)\subset \Z^++s$
for suitable $s$ than $A$ has at most $3$ elements. This result is motivated by the
conjecture that the disk has not more than $3$ orthogonal exponentials.

Further, we prove that if $A$ is a set of exponentials mutually orthogonal with respect to any
symmetric convex set $K$ in the plane with a smooth boundary and everywhere non-vanishing curvature,
then $\# (A \cap {[-q,q]}^2) \leq C(K) q$ where $C(K)$ is a constant depending only on $K$.
This extends and clarifies in the plane the result of Iosevich and Rudnev.
As a corollary, we obtain the result from \cite{IKP01} and \cite{IKT01} that if $K$ is a centrally symmetric
convex body with a smooth boundary and non-vanishing curvature, then $L^2(K)$ does not possess an orthogonal
basis of exponentials.
\end{abstract}

\maketitle

\tableofcontents

\section{Introduction}

The aim of this paper is to investigate size properties of a set $A$ whose distance set 
$$
\Delta(A):=\{|a-a'|\,: a,a'\in A\}
$$
has some prescribed arithmetic properties. For instance, Solymosi (\cite{S04}, {\it see also} \cite{S04b})
proved the following: 

\begin{theorem}[Solymosi, \cite{S04}]
\label{th:solymosi}\ \\
Let $A \subset \R^2$ such that $\Delta(A) \subset \Z^{+}$. Then 
\begin{equation}
\label{eq:0.6}
\# (A \cap {[-q,q]}^2) \lesssim q.
\end{equation}
\end{theorem}

This result is essentially sharp as can be seen in the following way. Let $A_N$ be the subset of the plane
consisting of $(n,0)$, where $n$ is a large positive integer, and pairs of the form $(0,m)$ such that $m$ is
a positive integer and $n^2+m^2=l^2$ for some positive integer $l$.
By elementary number theory, one can find approximately 
$\frac{N}{\sqrt{\log(N)}}$ such integers $m$ that are less than $N$. 

Our main result can be seen as a shifted version of this theorem:

\begin{theorem}\label{th:sols}\ \\
Let $A\subset\R^2$ be a set such that $\Delta(A)\subset\Z^++s$
where $s$ is either a trancendental number or a rational that is
not a half-integer nor an integer. Then $A$ has at most $3$ elements.
\end{theorem}

Our original motivation to the above question comes from questions related to the Fuglede
conjecture\footnote{Strictly speaking, this is no longer a conjecture as it has been disproved in high enough dimension by Tao \cite{Tao}, Matolcsi \cite{Ma} and
Kolountzakis-Matolcsi \cite{KMa}, see also \cite{KMb}} on the existence of orthonormal bases for given sets. To be more precise,
if $K$ is a compact set in $\R^d$, we will say that the collection $\{e^{2 \pi i \scal{x, a}}\}_{a \in A}$
is orthogonal with respect to $K$ if
\begin{equation}
\label{eq:0.1}
\int_K e^{2 \pi i x \cdot (a-a')} dx=0,
\end{equation}
whenever $a \not=a' \in A$. In other words, we ask for a set $A\subset\R^d$ such that
the difference set $A-A:=\{a-a'\,: a,a'\in A\}$ is included in in the zero set of
the Fourier transform $\widehat{\chi_K}$ of the characteristic function $\chi_K$ of $K$.

In the particular case where $K$ is the unit disc in the plane, $\widehat{\chi_K}$ is a radial function
and its zero set is the set of circles with radius the zeroes of the Bessel function $J_1$.
Further, up to re-scaling, it is known that these zeroes are of the form $k+\frac{1}{4}+error$ with $k\in\Z_+$.
Dropping the error term, we thus see that
Theorem \ref{th:sols} brings ground on a conjecture of Fuglede \cite{Fu}
that the disc has at most 3 orthogonal exponentials.

In the case of more general convex sets, many results are known. For instance, assume that
$K$ is a compact convex set in $\R^d$, symmetric with respect to the origin, such that the boundary
of $K$ is smooth and has everywhere non-vanishing curvature. It was proved in \cite{IKT01} that such sets
do not admit an orthonormal \emph{basis} of exponentials. Further, in \cite{IR03}, the authors proved that
$\{e^{2 \pi i \scal{x, a}}\}_{a \in A}$ is an infinite set of orthogonal exponentials,
then $A$ is contained in a line. Moreover, they
proved that if $d \not=1 \mod(4)$, then $A$ is necessarily finite. When $d=1 \mod(4)$ the authors
produced examples of convex bodies $K$ for which there exists infinite $A$'s satisfying (\ref{eq:0.1}).
Unfortunately, the proof in \cite{IR03} does not provide a finite upper bound for the size of $A$
leaving open the possibility of having arbitrarily large sets of orthogonal exponentials.

The reason for this is reliance on an asymptotic generalization of the following combinatorial principle
due to Anning and Erd\"os (\cite{Er45a} {\it see also} \cite{Er45}):

\begin{theorem}[Anning and Erd\"os, \cite{Er45a,Er45}]
\label{th:erdos}\ \\
Let $d\geq 2$, and let $|\cdot|$ be the Euclidean norm on $\R^d$.
For $A \subset \R^d$, let 
\begin{equation}
\label{eq:0.2}
\Delta(A)=\{|a-a'|: a,a' \in A\}.
\end{equation}
If $\# A=\infty$ and $\Delta(A) \subset \Z^{+}$, then $A$ is a subset of a line. 
\end{theorem}

The reason explicit bounds are difficult to extract from any application of this principle is that
for any $N$ there exists $A_N \subset \R^d$ not contained in a line such that $\# A_N=N$ and
$\Delta(A_N) \subset \Z^{+}$.  This suggests that a different geometric point of view is needed
to extract an explicit numerical bound if one exists. Moreover, such a point of view is likely to be
dimension specific because, as we mention, above, when $d=1 \mod(4)$, the set of orthogonal exponentials
may be infinite. 

Even though an explicit numerical upper bound still eludes us, we have been able to prove the following:

\begin{theorem}
\label{th:2}\ \\
Let $K$ be a convex planar set, symmetric with respect to the origin. Suppose that the boundary of
$K$ is smooth and has everywhere non-vanishing curvature. Then there exists a constant $C$
depending only on $K$ such that, whenever $A \subset \R^d$ is a set such that (\ref{eq:0.1}) holds, then
\begin{equation}
\label{eq:0.4}
\# (A \cap {[-q,q]}^d) \leq Cq^{d-1}.
\end{equation}
In particular, $K$ does not have an orthogonal basis of exponentials.
\end{theorem}

The constant $C(K)$ actually only depends on the inner and outer radius of $K$ as well as on bounds on its curvature.
The second part of the theorem is well known and was first proved in \cite{IKT01} (\cite{IKP01,Fu2} when $K$ is the Euclidean ball)
and also follows immediately
from \cite{IR03}. Our proof is not substantially different from previous proofs in the field. Nevertheless
it slightly simplifies the argument in \cite{IKT01} and allows to obtain a quantitative bound.
 
\medskip

This paper is organized as follows. We start by proving Theorem \ref{th:sols} and some generalizations of it.
We then complete the paper by devoting Section \ref{sec:1} to the proof of
Theorem \ref{th:2}.

\section{A shifted Erd\"os-Solymosi theorem}
In the remaining of the paper, we will identify $\R^2$ and $\C$.

\subsection{Proof of Theorem \ref{th:sols}}\ \\
In this section, we will prove the following slightly stronger
version of Theorem \ref{th:sols}

\begin{theorem}
\label{th:solymosis}\ \\
Let $s\in\R$ be such that $8s^8\notin \Z+4s\Z+2s^2\Z+4s^3\Z+s^4\Z+2s^5\Z+2s^6\Z+4s^7\Z$.
If $A \subset \R^2$ is such that $\Delta(A) \subset \Z^{+}+s$, then $\#A\leq 3$.
\end{theorem}

The result is best possible as if $a_0=0$, $a_1=r$, $a_2=re^{i\pi/3}$, then all distances are equal to $r$.

\begin{proof} 

Assume that $a_0,a_1,a_2,a_3$ are four different elements from $A$.
For $j=1,\ldots,3$ let $\alpha_j=a_j-a_0$ and write $\alpha_j=(k_j+s)e^{i\theta_j}$.
There is no loss of generality in assuming that $\alpha_1$ is real,
that is $\theta_1=0$.

For $j=2,3$ let $\beta_j=a_j-a_1=\alpha_j-\alpha_1$. Write $|\beta_j|=l_j+s$ with $l_j\in\Z^+$.

From $|\beta_j|^2=|\alpha_j-\alpha_1|^2=|\alpha_1|^2+|\alpha_j|^2-2|\alpha_1||\alpha_j|\cos\theta_j$, we deduce that
\begin{equation}
\label{eq:costhetaj}
2|\alpha_1||\alpha_j|\cos\theta_j=|\alpha_1|^2+|\alpha_j|^2-|\beta_j|^2.
\end{equation}
From this, we get that
\begin{equation}
\label{eq:sinthetaj}
(2|\alpha_1||\alpha_j|\sin\theta_j)^2=
4|\alpha_1|^2|\alpha_j|^2-(|\alpha_1|^2+|\alpha_j|^2-|\beta_j|^2)^2
\end{equation}

On the other hand $\alpha_3-\alpha_2=a_3-a_2$ thus $|\alpha_3-\alpha_2|=m+s$ for some $m\in\Z^+$ and
\begin{equation}
\label{eq:thetaj}
2|\alpha_2||\alpha_3|\cos(\theta_3-\theta_2)=|\alpha_2|^2+|\alpha_3|^2-|\alpha_3-\alpha_2|^2,
\end{equation}
from which we get that
\begin{eqnarray}
4|\alpha_1|^2|\alpha_2||\alpha_3|\cos(\theta_3-\theta_2)&=&2|\alpha_1|^2(|\alpha_2|^2+|\alpha_3|^2-|\alpha_3-\alpha_2|^2)\nonumber\\
&=&2(k_1+s)^2\bigl((k_2+s)^2+(k_3+s)^2-(m+s)^2)\label{ct23}\\
&=&2k_1^2(k_2^2+k_3^2-m^2)+4s\bigl(k_1^2(k_2+k_3-m)+k_1(k_2^2+k_3^2-m^2)\bigr)\nonumber\\
&&\qquad+s^2\bigl(k_1^2+k_2^2+k_3^2-m^2+4k_1(k_2+k_3-m)\bigr)\nonumber\\
&&\qquad\qquad+4s^3(k_1+k_2+k_3-m)+2s^4\nonumber\\
&\in& 2\Z+s4\Z+s^22\Z+s^34\Z+2s^4.\nonumber
\end{eqnarray}

On the other hand,
$$
4|\alpha_1|^2|\alpha_2||\alpha_3|\cos(\theta_3-\theta_2)=
4|\alpha_1|^2|\alpha_2||\alpha_3|(\cos\theta_2\cos\theta_3+\sin\theta_2\sin\theta_3).
$$
But, with (\ref{eq:costhetaj}),
\begin{eqnarray}
4|\alpha_1|^2|\alpha_2||\alpha_3|\cos\theta_2\cos\theta_3&=&\bigl((k_1+s)^2+(k_2+s)^2-(l_2+s)^2\bigr)\nonumber\\
&&\qquad\times\bigl((k_1+s)^2+(k_3+s)^2-(l_3+s)^2\bigr)\label{ct2t3}\\
&\in& \Z+s2\Z+s^2\Z+s^32\Z+s^4,\nonumber
\end{eqnarray}
%&=&(k_1^2+k_2^2-l_2^2)(k_1^2+k_3^2-l_3^2)\nonumber\\
%&&\quad+2s\bigl((k_1^2+k_2^2-l_2^2)(k_1+k_3-l_3)+(k_1^2+k_3^2-l_3^2)
%(k_1+k_2-l_2)\bigr)\nonumber\\
%&&\quad+s^2\bigl(2k_1^2+k_2^2+k_3^2-l_2^2-l_3^2+4(k_1+k_2-l_2)
%(k_1+k_3-l_3%)\bigr)\nonumber\\
%&&\quad+2s^3(2k_1+k_2+k_3-l_2-l_3)+s^4\nonumber\\
%&\in& \Z+s2\Z+s^2\Z+s^32\Z+s^4.\nonumber
%\end{eqnarray}
as can be seen by expanding the expression in the first line.
It follows that
$$
4|\alpha_1|^2|\alpha_2||\alpha_3|\sin\theta_2\sin\theta_3=
4|\alpha_1|^2|\alpha_2||\alpha_3|\cos(\theta_3-\theta_2)-
4|\alpha_1|^2|\alpha_2||\alpha_3|\cos\theta_2\cos\theta_3
$$
has to be in $\Z+s2\Z+s^2\Z+s^32\Z+s^4$, thus
$$
(4|\alpha_1|^2|\alpha_2||\alpha_3|\sin\theta_2\sin\theta_3)^2\in
\Z+s4\Z+s^22\Z+s^34\Z+s^4\Z+s^54\Z+s^62\Z+s^74\Z+s^8.
$$
Now, from (\ref{eq:sinthetaj}), we get that
\begin{eqnarray}
(2|\alpha_1||\alpha_j|\sin\theta_j)^2&=&
4(k_1+s)^2(k_j+s)^2-\bigl((k_1+s)^2+(k_j+s)^2-(l_j+s)^2\bigr)^2\label{sint2t3}\\
%&=&\bigl((k_1+k_2)^2-l_2^2\bigr)\bigl(l_2^2-(k_1-k_2)^2\bigr)\nonumber\\
%&&+4s\bigl(k_1^2(k_2+l_2-k_1)+k_2^2(l_2+k_1-k_2)+l_2^2(k_1+k_2-l_2)
%\bigr)\nonumber\\
%&&+2s^2(-k_1^2-k_j^2-l_j^2+4k_1k_2+4k_1l_2+4k_2l_2)\nonumber\\
%&&+4s^3(k_1+k_j+l_j)+3s^4\nonumber\\
&\in& \Z+s4\Z+s^22\Z+s^34\Z+3s^4\nonumber
\end{eqnarray}
as previously, thus
$$
(4|\alpha_1|^2|\alpha_2||\alpha_3|\sin\theta_2\sin\theta_3)^2\in
\Z+4s\Z+s^22\Z+s^34\Z+s^4\Z+s^54\Z+s^62\Z+s^712\Z+9s^8.
$$
We thus want that
$$
8s^8\in\Z+4s\Z+2s^2\Z+4s^3\Z+s^4\Z+2s^5\Z+2s^6\Z+4s^7\Z
$$
which contradicts our assumption on $s$.
\end{proof}

\begin{remark}
The assumption on $s$ is quite mild as it is satisfied by all transcendental numbers, all algebraic numbers of order at least 9
and also by
all rational numbers that are not integers nor half-integers. For this last fact, if $s=\frac{p}{q}$
with $p,q$ mutually prime and $q\not=1$, then the assumption reads
$$
8p^8\notin q^8\Z+4pq^7\Z+2p^2q^6\Z+4p^3q^5\Z+p^4q^4\Z+2p^5q^3\Z+2p^6q^2\Z+4p^7q\Z\subset q\Z
$$
so that $q$ divides $8$. But, then writing $q=2r$ with $r=1,2$ or $3$, the assumption reads
$$
8p^8\notin 2^8r^8\Z+2^9pr^7\Z+2^7p^2r^6\Z+2^6p^3r^5\Z+2^4p^4r^4\Z+2^4p^5r^3\Z+2^3p^6r^2\Z+2^3p^7r\Z\subset 8r\Z
$$
so that $r=1$ and $s=\frac{p}{2}$.

Also, note that one may scale the assumption. For example, if we assume that $\Delta(A)\subset \frac{1}{2}\Z^++\frac{1}{8}$
then $\#A\leq 3$, since $B=2A=\{2a,\ a\in A\}$ satisifes $\Delta(B)\subset\Z+\frac{1}{4}$, which establishes the link
with zeroes of the Bessel function $J_1$ ({\it see} next section).

Of course, we may scale both ways, and we then for instance get
that if $\Delta(A)\subset 4\Z^++1$ or if $\Delta(A)\subset 4\Z^++3$
(that is $s=1/4$ and $s=3/4$ respectively)
then $\#A\leq 3$. This result is false when $s=1/2$
since it is not hard to construct
sets for which  $\Delta(A)\subset 2\Z^++1$.
Indeed, let $k,l\in\Z^+$ and assume that $2l+1\leq 2(2k+1)$
and let $\theta=\arccos\frac{2l+1}{2(2k+1)}$.
Finally, let $a_0=0$, $a_1=1$, $a_2=(2k+1)e^{i\theta}$ and $a_3=-(2k+1)e^{-i\theta}$.
It is then clear that $a_i-a_0$, $i=1,2,3$ all have odd integer modulus.
Further $a_3-a_2=-2(2k+1)\cos\theta=-2l-1$ is an odd integer.
Finally
\begin{eqnarray*}
|a_2-a_1|
&=&\bigl((2k+1)\cos\theta-1\bigr)^2+(2k+1)^2\sin^2\theta\\
&=&(2k+1)^2+1-2(2k+1)\cos\theta=4k^2+2(2k-l)+1
\end{eqnarray*}
while
\begin{eqnarray*}
|a_3-a_1|
&=&\bigl(-(2k+1)\cos\theta-1\bigr)^2+(2k+1)^2\sin^2\theta\\
&=&(2k+1)^2+1+2(2k+1)\cos\theta=4k^2+2(2k+l)+3
\end{eqnarray*}
are also an odd integers. It seems nevertheless clear from the previous proof that some sparcity should happen in this case.
% as there relations in the difference set.
\end{remark}

\subsection{A perturbation of Theorem \ref{th:sols}}\ \\
Recall that our original motivation in proving Theorem \ref{th:sols} was to bound the number of
exponentials that are orthogonal for the disc in the plane. Recall the well known fact that
$\widehat{\chi_B}(\xi)=\frac{1}{|\xi|}J_1(2\pi|\xi|)$ where $J_1$ is the Bessel function
of order $1$. It immediately follows that $\{e^{2i\pi\scal{a,x}}\}_{a\in A}$
is an orthogonal set of exponentials for the disc if and only if the distance set
$\Delta(A)$ of $A$ satisfies $\Delta(A)\subset\mathcal{Z}_{J_1}$ where $\mathcal{Z}_{J_1}$
is the set of zeroes of $J_1$.

Further, $J_1$ has the following asymptotic expansion when $r\to+\infty$ ({\it see e.g.}
\cite[VIII 5.2, page 356-357]{St93}):
$$
J_1(r)\sim-\left(\frac{2}{\pi r}\right)^{1/2}\ent{
\sum_{j=0}^\infty\frac{(-1)^j\Gamma\left(\frac{3}{2}+2j\right)}
{2^{2j}(2j)!\Gamma\left(\frac{3}{2}-2j\right)}\frac{\sin\left(r-\frac{\pi}{4}\right)}{r^{2j}}
-
\sum_{j=0}^\infty\frac{(-1)^j\Gamma\left(\frac{5}{2}+2j\right)}
{2^{2j+1}(2j+1)!\Gamma\left(\frac{1}{2}-2j\right)}\frac{\cos\left(r-\frac{\pi}{4}\right)}{r^{2j+1}}}.
$$
From this, it is not hard to see that if $\widehat{\chi_B}(\xi)=0$ and if $\xi$ is big enough, then
\begin{equation}
\label{eq:asbeord0}
4|\xi|^2=(k+1/4)^2+\frac{3}{4\pi^2}+O(k^{-2}).
\end{equation}

We will now show that we can still perturbate Theorem \ref{th:sols}
so as that the difference set $\Delta(A)$ consists
of a small perturbation (and a harmless re-scaling) of the zeroes of the Bessel function in the following sense:

\begin{proposition}\ \\
Let $s\in\Q\setminus\frac{1}{2}\Z$ and let $\eta$ be 
either algebraic of order at least $5$ or trancendental.
Let $A\subset\R^2$ be such that every element $\alpha\in\Delta(A)$
has the property that $|\alpha|^2=(k+s)^2+\eta$,
then $A$ has at most $3$ elements.
\end{proposition}

\begin{proof} We use the same notation and assumptions as in the proof
of Theorem \ref{th:sols}.
Identity (\ref{ct23}) then becomes
\begin{eqnarray*}
4|\alpha_1|^2|\alpha_2||\alpha_3|\cos(\theta_3-\theta_2)&=&
2\bigl((k_1+s)^2+\eta\bigr)\bigl((k_2+s)^2+(k_3+s)^2-(m+s)^2+\eta)\\
&=&2k_1^2(k_2^2+k_3^2-m^2)+4s\bigl(k_1^2(k_2+k_3-m)+k_1(k_2^2+k_3^2-m^2)\bigr)\\
&&+s^2\bigl(k_1^2+k_2^2+k_3^2-m^2+4k_1(k_2+k_3-m)\bigr)\\
&&+4s^3(k_1+k_2+k_3-m)+2s^4\\
&&+2\eta\bigl((k_1+s)^2+(k_2+s)^2+(k_3+s)^2-(m+s)^2\bigr)+2\eta^2\\
&\in&\Q+\eta\Q+2\eta^2.
\end{eqnarray*}

Identity (\ref{ct2t3}) becomes
\begin{eqnarray*}
4|\alpha_1|^2|\alpha_2||\alpha_3|\cos\theta_2\cos\theta_3&=&
\bigl((k_1+s)^2+(k_2+s)^2-(l_2+s)^2+\eta\bigr)\\
&&\times\bigl((k_1+s)^2+(k_3+s)^2-(l_3+s)^2+\eta\bigr)\\
%&=&(k_1^2+k_2^2-l_2^2)(k_1^2+k_3^2-l_3^2)\\
%&&+2s\bigl((k_1^2+k_2^2-l_2^2)(k_1+k_3-l_3)
%+(k_1^2+k_3^2-l_3^2)(k_1+k_2-l_2)\bigr)\\
%&&+s^2\bigl(2k_1^2+k_2^2+k_3^2-l_2^2-l_3^2
%+4(k_1+k_2-l_2)(k_1+k_3-l_3)\bigr)\\
%&&+2s^3(2k_1+k_2+k_3-l_2-l_3)+s^4\\
%&&+\eta\bigl(2(k_1+s)^2+(k_2+s)^2+(k_3+s)^2
%-(l_2+s)^2-(l_3+s)^2\bigr)+\eta^2\\
&\in&\Q+\eta\Q+\eta^2.
\end{eqnarray*}
It follows that $4|\alpha_1|^2|\alpha_2||\alpha_3|\bigl(\cos(\theta_3-\theta_2)-\cos\theta_2\cos\theta_3\bigr)$
belongs to
$$
\Q+\eta\Q+\eta^2.
$$
Squaring, we get that $\bigl(4|\alpha_1|^2|\alpha_2||\alpha_3|\sin\theta_2\sin\theta_3\bigr)^2$ is in
$$
\Q+\eta\Q
+\eta^2\Q+\eta^3\Q+\eta^4.
$$

On the other hand, Identity (\ref{sint2t3}) becomes
\begin{eqnarray*}
(2|\alpha_1||\alpha_j|\sin\theta_j)^2&=&
4\bigl((k_1+s)^2+\eta\bigr)\bigl((k_j+s)^2+\eta\bigr)\\
&&-\bigl((k_1+s)^2+(k_j+s)^2-(l_j+s)^2+\eta\bigr)^2\\
%&=&\bigl((k_1+k_j)^2-l_j^2\bigr)\bigl(l_j^2-(k_1-k_j)^2\bigr)\\
%&&+4s\bigl(k_1^2(k_j+l_j-k_1)+k_2^2(l_j+k_1-k_j)
%+l_j^2(k_1+k_j-l_j)\bigr)\\
%&&+2s^2(-k_1^2-k_j^2-l_j^2+4k_1k_j+4k_1l_j+4k_jl_j)\\
%&&+4s^3(k_1+k_j+l_j)+3s^4\\
%&&+\eta\bigl(2(k_1+s)^2+2(k_j+s)^2-2(l_j+s)^2\bigr)+3\eta^2\\
&\in& \Q+\eta\Q+3\eta^2.
\end{eqnarray*}
It follows that $\bigl(4|\alpha_1|^2|\alpha_2||\alpha_3|\sin\theta_2\sin\theta_3\bigr)^2$ is in
$\Q+\eta\Q+\eta^2\Q+\eta^3\Q+9\eta^4$.
As we have assumed that $\eta\Q+\eta^2\Q+\eta^3\Q+\eta^4\Q\subset\R\setminus\Q$,
we have obtained the contradiction $8s^8\in\Q+\R\setminus\Q\subset\R\setminus\Q$.
\end{proof}

\begin{remark}
We have used that the perturbation by $\eta$ is fixed only in a mild way
in order to simplify computations.
Actually, it is not hard to see that if we assume that
if each $\alpha\in\Delta(A)$ has the property that
$|\alpha|^2=(k+s)^2+P_k(\eta)$ where $k$ is an integer,
$\eta$ is a fixed \emph{trancendental} number and $P_k$
is a polynomial with rational coefficients, then the above proof
still gives the same result. If moreover the degrees of the
$P_k's$ are bounded by $M$, then the proof stille works provided
$\eta$ is algebraic of order at least $4M+1$.

In particular, let us recall that the asymptotic expansion \eqref{eq:asbeord0}
of the large zeroes $\xi$ of $\widehat{\chi_B}$ can be pushed further
to obtain:
\begin{equation}
|\xi|^2=\frac{1}{4}(k+1/4)^2\left(1+\frac{3}{4\pi^2(k+1/4)^2}+\sum_{j=2}^N\frac{c_j}{\pi^{2j}(k+1/4)^{2j}}+O\left(\frac{1}{k^{2N+2}}\right)\right).
\label{eq:asympt}
\end{equation}
where the $c_j$'s are rational constants. Let us now truncate this formula. More precisely,
let us assume that the set $A$ is such that each $\alpha\in\Delta(A)$ has the property that $|\alpha|^2$
is of the form
$$
\frac{1}{4}(k+1/4)^2
\left(1+\sum_{j=1}^N\frac{c_j}{\pi^{2j}(k+1/4)^{2j}}\right).
$$
Then $A$ has at most $3$ elements.

An even more careful examination of the proof shows that,
if each $\alpha\in\Delta(A)$ is of the form $|\alpha|^2=(k+s)^2+\eta_k$
then only $6$ $\eta_k$'s intervene in the proof
(corresponding to $k_1,\ldots,k_3$, $l_2,l_3$ and $m$).
Moreover they are raised to the power at most $4$, so $A$ has at most $3$ elements as soon as there
exists no rational polynomial relation of degree at most $4$
between any $6$ $\eta_k$'s {\it i.e}
if, for any $j_1,\ldots,j_6$, the only polynomial of degree $\leq 4$ of $6$ variables
with coefficients in $\Q$ such that $P(\eta_{j_1},\ldots,\eta_{j_6})=0$
is $P=0$. Such relations are highly unlikely between zeroes of the
Bessel function $J_1$, it is thus natural to conjecture, following Fuglede \cite{Fu} that the disk has no more than $3$ orthogonal exponentials.
\end{remark}

\section{Orthogonal exponentials for planar convex sets}
\label{sec:1}

For sake of simplicity, we will concentrate on the proof of Theorem \ref{th:2} in the case of dimension 2.

\subsection{Preliminaries}\ \\
We will need the following well known facts about convex sets.

\begin{notation}
For a convex set $K$, we call $\rho_K$ its Minkowski function of $K$, so that $K=\{x: \rho(x) \leq 1\}$,
and $\rho_K^{*}$ its support function given by
\begin{equation}
\rho_K^*(\xi)=\sup_{x \in K} \scal{x,\xi}.
\label{eq:1.3}
\end{equation} 
\end{notation}

By the method of stationary phase (\cite{He62}, {\it see e.g.} \cite[Chapter 3]{St93}), 
\begin{equation}
\widehat{\chi}_K(\xi)=C_1{|\xi|}^{-\frac{3}{2}} \sin \left(2 \pi \left(\rho_K^{*}(\xi)-\frac{1}{8} \right) \right)
+E(\xi),
\label{eq:1.1}
\end{equation}
with 
\begin{equation}
|E(\xi)| \leq C_2{|\xi|}^{-\frac{5}{2}},
\label{eq:1.2}
\end{equation}
where $C_1$ and $C_2$ are some constants depending only on $K$.

It should also be noted that if $\{e^{2i\pi ax}\}_{a\in A}$ is 
orthogonal with respect to $L^2(K)$ then $\widehat{\chi_K}(a-a')=0$ for 
$a,a'\in A$. But $\widehat{\chi_K}$ is continuous and 
$\widehat{\chi_K}\not=0$ so there exists $\eta_0$ such that $|a-a'|\geq 
\eta_0$, that is, the set $A$ is separated with separation depending 
only on $K$.

\subsection{Proof of Theorem \ref{th:2}}\ \\
An immediate consequence of \eqref{eq:1.1} and \eqref{eq:1.2}
is that, if $A$ is as in the statement of Theorem \ref{th:2},
there exists a constant $C_3$ such that, whenever $a,a' \in A$, then 
\begin{equation}
\left|\rho^{*}_K(a-a')-\frac{k}{2}-\frac{1}{8}\right|\leq\frac{C_3}{k+1}
\label{eq:1.4a}
\end{equation} 
for some integer $k$. We may now cut $A$ into a finite number of pieces,
such that in each piece, any two elements $a,a'$
are separated enough to have $k\geq 100C_3$ in \eqref{eq:1.4a}.

Now, if $a,a',a''\in A$ are in a $q$ by $\alpha$ rectangle $\mathcal{R}_\alpha$
then $a-a',a-a'',a'-a''$ are all in an angular sector
with direction some vector $e$ that depends only on
$\mathcal{R}_\alpha$. More precisely, the angle $\theta=\theta(e,a-a')$
between $e$ and $(a-a')$ is at most $\theta$ with
$\sin\theta=\alpha/|a-a'|$. In particular,
$\theta\leq\theta_m$  where $\sin\theta_m=\alpha/L$ and $L$ is the 
minimal distance between two elements of $A$. Further, from the 
curvature assumption
on $K$, for $u$ in such a sector,
$$
\bigl|\rho_K^*(u)-c|u|\bigr|\leq c'\theta(u,e)^2|u|,
$$
where $c=c(e)$ and $c'=c'(e)$ are two consants that depend on $e$,
provided $\theta(u,e)$ is small enough (that is $\alpha$ is taken to be small enought). In particular, for $u=a-a'$,
$$
\bigl|\rho_K^*(a-a')-c|a-a'|\bigr|\leq \frac{C_4}{|a-a'|}.
$$
It follows from this and \eqref{eq:1.4a}, that there exists an integer
$k$ such that
\begin{equation}
\left|c|a-a'|-\frac{k}{2}-\frac{1}{8}\right|\leq\frac{C_5}{k+1}.
\label{eq:1.4b}
\end{equation} 
There is no loss of generality to assume that elements in $A$
are sufficiently separated to have $C_5/(k+1)<1/100$.
Similarely, there also exist integers $l,m$ such that
\begin{equation}
\abs{c|a-a''|-\frac{l}{2}-\frac{1}{8}}\leq \frac{1}{100}
\quad\mbox{and}\quad
\abs{c|a'-a''|-\frac{m}{2}-\frac{1}{8}}\leq \frac{1}{100}.
\label{eq:1.4c}
\end{equation}

Now, since $a,a',a''$ are in a box of size $q$ by $\alpha$,
if $|a-a''|\geq |a-a'|,|a'-a''|$ then
\begin{eqnarray*}
|a-a''|&\geq&(|a-a'|^2-\alpha^2)^{1/2}+ (|a'-a''|^2-\alpha^2)^{1/2}\\
&=&|a-a'|+|a'-a''|-\frac{\alpha^2}{(|a-a'|^2-\alpha^2)^{1/2}+|a-a'|}
-\frac{\alpha^2}{(|a'-a''|^2-\alpha^2)^{1/2}+|a'-a''|}\\
&\geq&|a-a'|+|a'-a''|-\frac{1}{100c}
\end{eqnarray*}
where $c$ is the constant in \eqref{eq:1.4b}, provided we have taken $\alpha$ small enough. It follows that
$$
\bigl|c|a-a''|-c|a-a'|-c|a'-a''|\bigr|\leq\frac{1}{100}.
$$
But the, from \eqref{eq:1.4b} and \eqref{eq:1.4c},
$$
\abs{\frac{l-m-k}{2}-\frac{1}{8}}\leq \frac{1}{25}
$$
a clear contradiction. Thus every $q$ by $\alpha$ rectangle
contains at most $2$ elements of $a$. It follows that
$A$ has at most $2q/\alpha$ elements in a $q\times q$ square.

\begin{remark}
The same proof in dimension $d\not=1\mod(4)$ works provided we use 
$q\times\alpha\times\cdots\times\alpha$ tubes. We would then obtain
that $A$ has at most $\lesssim q^{d-1}$ elements in any cube
of side $d$.
\end{remark}

\subsection{Orthogonal exponential bases}\label{sec:3}\ \\
The following result is proved in \cite{IKP01} in the case of the ball, and in
\cite{IKT01} in the general case. We shall give a completely self-contained and transparent proof below. 

\begin{theorem}[Iosevich, Katz, Pedersen, Tao, \cite{IKP01,IKT01}]
\label{th:3.1}\ \\
Let $K$ be a symmetric convex set in $\R^d$ with a smooth boundary and everywhere non-vanishing curvature.
Then $L^2(K)$ does not possess an orthogonal basis of exponentials.
\end{theorem} 

\begin{proof}
To prove Theorem \ref{th:3.1}, assume that $L^2(K)$ 
does possess an orthogonal basis of exponentials
${\{e^{2 \pi i x \cdot a}\}}_{a \in A}$.
From Theorem \ref{th:2}, $\# (A \cap [-q,q]^d) \lesssim q^{d-1}$.

But, as is well known \cite{Be66,La67,GR96,IKP01}, if ${\{e^{2 \pi i x \cdot a}\}}_{a \in A}$ is an orthonormal basis
of exponentials of $L^2(K)$, then $\dst\limsup\frac{\# (A \cap [-q,q]^d)}{q^d}>0$, a contradiction.
\end{proof}

\bibliographystyle{plain}

\end{document}